\documentclass[11pt]{article}

\usepackage[margin=1in]{geometry}
\usepackage[T1]{fontenc}
\usepackage[utf8]{inputenc}
\usepackage{lmodern}
\usepackage{microtype}
\usepackage{amsmath,amssymb}
\usepackage{array}
\usepackage{booktabs}
\usepackage{enumitem}
\usepackage{tabularx}
\usepackage{xcolor}
\usepackage[colorlinks=true,linkcolor=black,urlcolor=blue!60!black,citecolor=blue!60!black]{hyperref}

\setlist[itemize]{leftmargin=1.45em,topsep=3pt,itemsep=3pt}
\setlist[enumerate]{leftmargin=1.65em,topsep=3pt,itemsep=4pt}
\newcolumntype{Y}{>{\raggedright\arraybackslash}X}

\newcommand{\term}[1]{\emph{#1}}

\title{\textbf{Graduate Mathematics in the Age of AI}\\[0.35em]
\large Forming Mathematicians for Original, Independent, and Responsible Inquiry}
\author{Bacim Alali\thanks{bacimalali@math.ksu.edu}}

\date{}

\begin{document}

\maketitle

\begin{abstract}
Artificial intelligence can increasingly produce plausible, technically sophisticated
mathematical material faster than a developing graduate student can understand or verify
it. A sophisticated result or paper draft consequently becomes weaker evidence of the
student's own mathematical development. This creates a formation gap between output and
personal capacity, and a trust gap between a convincing argument and a warranted basis
for accepting it. The formation gap can persist even when the student genuinely
understands the output: understanding a supplied argument does not by itself establish
the capacity to initiate and direct inquiry.

These gaps are not the whole educational story. AI can also help students explore more
examples, compare approaches, enter unfamiliar areas, and undertake ambitious research.
The task is to design an apprenticeship that realizes these possibilities while
developing substantive mathematical command. The central purpose of a mathematics PhD
should be the formation of mathematicians capable of original, independent, and
responsible inquiry, including inquiry conducted in collaboration with AI.

This document develops that objective through four connected capacities: competence,
judgment, independence, and responsibility. It distinguishes a work's contribution to
mathematics from the evidence it provides of a student's formation; explains how a known
answer can initiate rather than end creative inquiry; and proposes changes in learning
activities, assessment, doctoral originality, advising, and institutional support.
Purposeful independent work and ambitious AI-assisted research are complementary parts
of the model. Its practical recommendations include proportionate contribution
statements, explicit recognition of advising costs, and staged pilots that evaluate
both internal mathematical ability and effective human--AI collaboration. The aim is
not to preserve an inherited sequence of training, but to improve mathematical formation
as the means of doing mathematics change.

\end{abstract}

\newpage
{\small
\tableofcontents
}
\newpage

\part{The apprenticeship under pressure}

\section{Purpose and premise}

What should graduate education in mathematics accomplish when AI can contribute
substantially to mathematical research? AI can assist with reading, computation,
conjecture formulation, proof construction, and exposition. These capabilities raise
questions that go beyond whether students should be permitted to use particular tools.
They require us to reconsider how graduate study develops mathematical understanding,
creativity, judgment, and independence---and how an apprenticeship should be designed
when producing a result, understanding it, and learning to initiate and direct research
need not develop together.

The discussion begins from a deliberately moderate premise:
\vspace{-0.1cm}
\begin{quote}
AI can increasingly produce plausible, technically sophisticated, research-level
mathematical material faster than a developing student can understand and check it.
\end{quote}
\vspace{-0.1cm}

This premise does not require the claim that present systems reliably solve arbitrary
research problems.  AI-generated proofs may contain hidden assumptions, fabricated
references, irrelevant detours, subtle gaps, or fatal errors.  That unreliability is part
of the educational problem rather than an objection to it.  Once a system can generate
material whose apparent sophistication exceeds a student's capacity to evaluate it, the
traditional connection between performance and formation is already under pressure.

The same capabilities create educational opportunities. Responsive explanation,
comparisons among approaches, computational exploration, and assistance with unfamiliar
prerequisites may help a student become a more capable mathematician. These benefits are
not automatic, any more than intellectual dependence is inevitable. They depend on what
the student does, what the system contributes, and how the activity is designed.

The aim is neither to preserve the current PhD unchanged nor to predict a single future.
It is to ask how graduate education can cultivate mathematical imagination, understanding,
technical power, and independence while making productive use of AI. The proposed
practices are designs to evaluate, not claims that one sequence of learning is best for
every student, field, or stage.

\section{The traditional apprenticeship}

The relationship between a mathematics PhD student and an advisor often develops through
two overlapping stages. Their balance matters because AI may change both the work of
formation and the research collaboration it makes possible.

\subsection{Formation before productivity}

Early in the apprenticeship, the student takes courses, learns an area, reads selected
papers, follows references, computes, constructs examples, reconstructs proofs, and
discusses difficulties with an advisor. The advisor identifies important gaps, chooses
worthwhile problems, and helps the student recognize what a successful argument would
require.

An assigned problem can be educationally successful even when it is not solved. An
attempt may reveal that an estimate fails at the boundary of its claimed parameter range,
that an analogy is superficial,
or a lemma has indispensable hypotheses. Failed proofs, partial arguments, and
well-formed questions can be evidence of developing understanding and research skill.
Few publishable results need not imply little progress.

Difficulty and time spent stuck are not educational achievements in themselves.
Trying an approach and discovering why it fails can build mathematical capacity;
remaining confused because of poor guidance, inaccessible explanations, or an
unrecognized prerequisite gap may not. A better apprenticeship should preserve
opportunities for the former while providing help with the latter. Human understanding requires intellectual activity and consolidation, but the
familiar duration and sequence of that activity are not beyond improvement.

\subsection{From student to contributor}

As students gain command of an area, they read independently, recognize proof
strategies, and make sustained contributions. They may identify a missing case,
introduce a technique, formulate a better question, or open a direction the advisor had
not anticipated. The relationship can become a collaboration and eventually support
an independent research program.

This creates a traditional apprenticeship bargain: an early investment in formation
can lead to later research collaboration, intellectual continuity, and a future
colleague. The relationship is not merely transactional. Advising serves an educational
mission even when a student does not advance the advisor's own project, and students
contribute through teaching, new perspectives, and the life of the department.

AI may alter this balance in more than one direction. It may reduce an advisor's reliance
on a student's routine research contributions, but it may also help a developing student
engage sooner in substantive collaboration. Neither the research benefits nor the
mentoring costs should be assumed in advance.

\section{The formation gap}

In the traditional apprenticeship, mathematical formation and research production were
imperfectly but substantially coupled. Working toward a result required students to
recognize gaps, read and experiment, attempt proofs, receive criticism, and revise their
ideas. The route to a result was also a route toward the capacity to understand and
extend it.

AI makes a shorter route to a mathematical artifact possible: a student can give an AI
system a research question and receive a candidate solution or polished draft before
undertaking comparable intellectual work. This creates a \term{formation gap}: the
distance between the sophistication of the available mathematical output and the
student's capacity to understand, evaluate, and develop mathematics, including
initiating and directing an investigation.

Consider a student who asks AI for help with a research problem and receives a useful
approach, a key lemma, and its proof. The student studies the argument carefully, fills
in details, and comes to understand it genuinely. When meeting with the advisor, the student
explains the proof convincingly and answers the advisor's questions. The project has
advanced, and real learning has occurred. Yet the student may not have practiced
formulating the lemma, choosing among uncertain approaches, recognizing why a path is
failing, or deciding what to try next. If this pattern recurs, those capacities may
receive little exercise while progress and convincing explanations suggest to the
advisor that the apprenticeship is developing them.

The concern is not that understanding supplied mathematics has no formative value.
Studying a proof can develop technique, imagination, and insight, and may itself require
substantial creative effort. The concern is that these gains do not automatically
establish the student's capacity to conduct an investigation before its direction or
answer has been supplied. \emph{The formation gap can persist after the understanding
gap has closed.} A good explanation is evidence of understanding, not sufficient evidence
of the full range of research capacities.

Presenting generated ideas as one's own adds an integrity problem, but removing that
misrepresentation does not remove the educational problem. The student may disclose
every substantive contribution, including assistance from AI, the proof may be correct,
and the explanation may be
excellent, while opportunities to formulate, search, choose, and revise remain repeatedly
bypassed. This is a risk to investigate, not a diagnosis of incapacity from AI use alone.

Two judgments must therefore remain distinct: \emph{what the work contributes to
mathematics}, and \emph{what the student's participation reveals about their formation}.
An important, reliably established theorem does not lose its mathematical value because
its proof was produced quickly or with AI. Its value, however, does not automatically
certify the maturity of the student presenting it. Conversely, an unsuccessful attempt
can develop substantial mathematical capacity without yielding a publishable result.

The appropriate response is not to restore difficulty wherever AI removes it or to
require rediscovery of every supplied idea. AI can reorganize the formative process:
a student might compare approaches, explore a richer class of examples, investigate a
connection, or start from a known answer and seek a new explanation or extension.
These activities become occasions for research formation when the student does
substantive mathematical work and learns to make and revise consequential choices.

The educational question is therefore not only whether the student understands the
result, but where reasoning, experimentation, choice, and invention occur across the
apprenticeship. Their distribution among student, advisor, and AI should be deliberately
designed for development, rather than left to whichever workflow produces answers most
quickly.

\section{The trust gap}

Suppose a student encounters a difficult passage in a reputable paper selected by an
advisor.  Although the paper may contain mathematical errors---refereeing is not
infallible---the student ordinarily has a reasonable working presumption that the passage
is mathematically sound.  If the student does not initially understand it, sustained
effort is likely to be worthwhile.  The presumption of soundness encourages the student
to treat confusion as something to work through, thereby turning it into productive
mathematical effort.

When the passage is generated by AI, the student may face at least four possibilities:

\begin{enumerate}
    \item the argument is correct but exceeds the student's present understanding;
    \item the argument is basically sound but contains a repairable gap;
    \item the argument is false and should be rejected or substantially rebuilt;
    \item the argument is formally sophisticated but irrelevant, vacuous, or
    mathematically uninformative.
\end{enumerate}

The student may lack exactly the knowledge needed to tell these cases apart.  Continuing
to study a false argument can waste time and create misconceptions; abandoning a correct
but difficult argument can prevent learning.  Asking AI for clarification may produce
another confident explanation rather than independent evidence.

Mathematical trust was never automatic.  It was produced by identifiable authorship,
editorial selection, refereeing, seminars, citation, correction, independent checking,
and the gradual evaluation of published work by a community.  These practices did not
guarantee correctness, but they gave readers reasons for provisional trust and routes for
challenging a claim.

AI-generated mathematics can reach a student before it has passed through comparable
checks. Existing vetting practices were designed for public mathematical claims, not for
the volume of individualized, real-time output produced during an AI interaction. The
student may therefore lack a reliable initial basis for deciding whether the difficulty
lies in their understanding or in the argument. The distinction is between evidence
and insufficient evidence, not simply between human and machine authorship: a carefully
checked AI-assisted proof may deserve more trust than an unchecked human argument.

A confidence estimate from another AI system may help if it has been calibrated against
known examples containing relevant kinds of errors.  A confidence label alone, however,
cannot create warranted trust.  A sound learning environment combines reliable
references; provenance, meaning a record of where a claim came from and how it was
transformed; independent derivations; tests of special cases; counterexamples; symbolic
or numerical checks; formal verification where available; peer criticism; and expert
review.  Students must also calibrate their own confidence, distinguishing between not
having found an error and having strong evidence that an argument is sound.  The
construction and evaluation of mathematical trust must become an explicit part of
graduate training.

The kinds of assurance must also be distinguished. Agreement among AI models need not be
independent confirmation; examples and numerical tests can expose errors without proving
a general claim. A formal certificate establishes a particular formal statement relative
to its assumptions and trusted checking system, while the correspondence between that
statement and the intended mathematics still requires scrutiny. Students need to learn
what each check establishes, not merely accumulate checks.

\section{The challenges for students, advisors, and graduate programs}

The formation and trust gaps create connected challenges at three levels. Students must
turn powerful assistance into their own mathematical development; advisors must design
and assess a formative apprenticeship; programs must sustain that apprenticeship and
justify what their degrees certify. The difficulties are not confined to academic
integrity. They arise even when AI use is fully disclosed and the resulting mathematics
is correct.

\subsection{For the student}

The student's challenge is to benefit from increasingly capable assistance without
becoming dependent on it for the direction and substance of every mathematical activity.
An explanation may feel clear while it is being read, yet leave the student unable to
use the idea in a different setting. More subtly, the student may understand it deeply
while repeatedly leaving the formulation of questions and the choice of approaches to
AI. The challenge is to develop both command of supplied mathematics and the ability
to begin and sustain inquiry when a promising route has not yet been given.

The student also faces decisions about when to seek help, which suggestions to pursue,
and how much of an unfamiliar argument to investigate. The trust gap complicates those
decisions: difficulty may reflect a missing prerequisite, a subtle mathematical idea,
or a false generated claim. At the same time, assistance can make unfamiliar areas and
ambitious projects more accessible. The challenge is to use that expanded opportunity
to develop technical skill, imagination, judgment, and increasing ownership of a
research direction, rather than merely the ability to obtain a finished answer.

\subsection{For the advisor}

AI changes the evidence available to the advisor, the design of student work, and
potentially the balance between research collaboration and educational investment.
Four challenges follow.

First, \emph{progress and genuine understanding are incomplete evidence of
research formation}. A student may explain an AI-supplied proof thoroughly while gaining
little experience in finding a direction before the answer is available. The advisor
must distinguish learning to understand an argument from learning to formulate
questions, choose approaches, and revise a plan under uncertainty. These capacities
overlap, but evidence of one cannot simply stand for all. The challenge is to observe
their development without turning every meeting into an examination or an investigation
of misconduct.

Second, \emph{the division of intellectual work must be designed for learning}.
AI may supply an immediate candidate solution to a problem selected to develop research
initiative. Reconstruction and explanation remain valuable, but may not exercise the
choices the assignment was intended to develop. The advisor must decide when the student
should attempt a formulation before receiving suggestions, when AI should extend or
criticize the student's approach, and when it is better to begin from an established
answer and investigate its explanation, consequences, or possible extensions.
This last possibility is \emph{answer-first inquiry}, developed in
Section~\ref{sec:answer-first}. Continually finding problems beyond AI's reach is neither stable
nor necessary. What matters is that students repeatedly practice initiating and
directing inquiry, with and without assistance, rather than only understanding work
whose direction has already been chosen for them.

Third, \emph{the material requiring expert attention may expand beyond the time available}.
A student may bring several generated approaches, lengthy proofs, or techniques unfamiliar
to both student and advisor. Assessing their reliability and educational value can require
substantial work. The advisor must decide what deserves close review, what the student
should investigate first, and when another specialist or verification tool is needed.
Without such decisions, advising risks becoming an unbounded review service. AI may
also reduce routine explanatory work or assist with checking; the net effect on
workload cannot simply be assumed.

Fourth, \emph{the practical incentives supporting intensive mentoring may change}.
If AI can perform some tasks previously undertaken by junior collaborators, an advisor
may become less dependent on those contributions while still needing to invest heavily
in a student's formation. One practical incentive for accepting that investment could
weaken. Conversely, AI may enable students to contribute more substantially and sooner,
making collaboration more valuable. The challenge is to sustain serious mentoring
without assuming either outcome. The educational value of advising does not depend on
a student's usefulness to the advisor's own research, but workload and reward structures
can still affect whether faculty have the capacity and incentive to undertake it.

Assigning students to check AI-generated output may appear to address both the review burden
and the need for useful research contributions. Verification can be intellectually
substantive and formative. If routine checking becomes the student's primary role,
however, the apprenticeship may develop neither question formulation nor ownership of a
research direction. This is why the response cannot consist simply of making students
quality-control workers for an advisor's AI-assisted research.

\subsection{For the program}

The program faces a challenge of both certification and institutional capacity. If
similar-looking dissertations can be produced through very different levels of student
understanding and participation, the artifact alone provides an inadequate basis for
certifying mathematical maturity. The program must distinguish a work's contribution
to mathematics from the evidence it provides of the student's capacity for independent
research. It must also decide what students should understand internally, what they may
do with assistance, and what constitutes justified reliance on others.

Those standards must be achievable under the conditions the institution actually
provides. More demanding assessment or mentoring cannot simply be added to existing
workloads without considering time, expertise, access, and cost. Unequal resources or
advising support can make the same formal requirements substantially different in
practice. The program's challenge is therefore to make the degree credible while
supporting a fair and sustainable route to earning it, with honest expectations about
duration and professional opportunities.

The later sections develop responses to these challenges: the redesign of advising in
Section~\ref{sec:advisor-model}, evidence of formation in
Section~\ref{sec:assessing-capacities}, and institutional recognition of mentoring in
Section~\ref{sec:advising-work}.

\clearpage
\part{Why human mathematical capacity still matters}

\section{The case for continued formation}

If AI can generate increasingly sophisticated mathematical results, proofs, and papers,
why should universities continue to invest years of faculty time and institutional
resources in developing human mathematicians?  Three complementary answers are needed.
The first concerns inquiry: trained mathematicians can contribute creatively to
discovery, judgment, and synthesis within an AI-assisted research environment.  The second is institutional and
normative: universities and mathematical communities must preserve the human capacity to
understand, teach, evaluate, and take responsibility for mathematical knowledge.  The
third concerns continuity and resilience: expertise takes years and several generations
to develop and transmit, so maintaining it preserves society's ability to respond to an
uncertain technological future.

\subsection{Discovery, judgment, and synthesis}

Deep mathematical training supports invention as well as evaluation. It helps researchers
formulate questions, recognize structure, devise representations, construct examples,
and develop methods. In an AI-assisted environment, these activities may be pursued with
a wider range of computational experiments and candidate approaches. The ambition is to
prepare mathematicians to use that expanded range intelligently and creatively.

Evaluation and synthesis also remain indispensable functions: someone or some system
must distinguish a significant connection from a superficial analogy, assess
explanatory power, and relate results to a coherent understanding of a field. If
AI-assisted production of candidate proofs and results expands faster than reliable
checking and incorporation into mathematical knowledge, these
functions become bottlenecks. It does not follow that their growing importance
guarantees a growing demand for human labor.

Graduate education should allow for the possibility that AI will contribute increasingly
to judgment, question selection, explanation, and synthesis as well as proof construction.
The extent and timing of such developments remain uncertain. The argument for graduate
formation should therefore not rest on a permanent division in which machines produce
and humans judge. It rests in part on the
value of people who can participate substantively in mathematical inquiry: form their
own questions, work with powerful assistance, assess competing reasons, and develop or
redirect a project.

The personal value of understanding and discovery matters too. Mathematics is a human
intellectual practice, not only a supply of useful results. Educating people to
participate deeply in it is a legitimate university purpose. This does not by itself
settle the appropriate scale or cost of graduate programs, but it is a reason that
cannot be reduced to a forecast of which tasks humans will perform more cheaply.

\subsection{Institutional stewardship and intellectual agency}

Human mathematical capacity also sustains functions that cannot be measured by paper
production alone.  Mathematicians referee results, teach the next generation, and
maintain standards concerning rigor, proof, evidence, attribution, and the responsible
communication of claims.  They evaluate the mathematical reliability of AI systems,
computational methods, software, and models; serve on scientific, professional, and
public panels; and detect when confidently presented errors begin to spread through a
research or professional community.

When mathematical claims influence engineering, science, finance, security, medicine, or
public decisions, qualified people must be able to determine what has been assumed, what
has been established, what remains uncertain, and what consequences may follow if the
mathematics is wrong.  A community that can request mathematical output but cannot
understand, criticize, or independently evaluate it has lost an important form of
intellectual agency.

Universities therefore have a stewardship responsibility: important mathematical
knowledge must remain open to human understanding, teaching, criticism, and
accountability.  Stewardship does not require humans to execute every calculation or
construct every proof without technological assistance.  It requires institutions to
retain sufficient expertise to understand consequential mathematical work, evaluate the
systems that produce it, challenge unreliable conclusions, and take responsibility for
decisions based on those conclusions.

\subsection{Continuity, resilience, and option value}

Mathematical expertise develops through years of study, research, teaching, and
collaboration. Once the communities that transmit it have substantially weakened,
preserved books and databases do not by themselves restore that capacity. This gives
continued formation an \term{option value}: it maintains the ability to respond as
future needs become clearer.

That capacity may matter when systems reveal unexpected limitations, access changes,
new scientific problems arise, or consequential applications require independent
assessment. It also preserves the ability to develop forms of human--AI collaboration
that cannot yet be anticipated. This is not a prediction of a permanent division of
labor. It is a commitment to the knowledge and institutional flexibility required to
adapt, alongside an obligation to examine the costs and outcomes of that commitment.

\section{A new educational objective}
\label{sec:capacities}

The central purpose of a mathematics PhD should be to develop the capacity for
\emph{original, independent, and responsible mathematical inquiry}. New results remain
important, but neither artifact production nor error detection alone defines the
objective. A mature mathematician can ask questions, develop ideas, bring a project to
fruition, explain its value, and establish a warranted basis for its claims.

Four connected capacities give this objective practical content:

\begin{enumerate}
    \item \textbf{Competence.} The student has substantive mathematical knowledge and
    technical ability. They can formulate precise statements, construct and analyze
    arguments, compute, build informative examples and counterexamples, and use
    definitions and methods in unfamiliar settings. Competence is productive as well
    as reconstructive: it supports the invention and development of mathematics.

    \item \textbf{Judgment.} The student distinguishes conjectures, heuristic arguments
    that suggest a conclusion without proving it, and rigorously established claims;
    evaluates assumptions and evidence; compares approaches;
    and assesses novelty, significance, explanatory power, and usefulness. They can
    give reasons for pursuing one question rather than another, while recognizing
    uncertainty and the possibility of revising that choice.

    \item \textbf{Independence.} The student identifies what they need to learn,
    formulates worthwhile questions, designs a research strategy, and selects among
    literature, collaborators, computation, AI, and formal tools. They can initiate,
    develop, and redirect inquiry rather than merely execute a supplied plan.
    Independence is compatible with substantial assistance and collaboration.

    \item \textbf{Responsibility.} The student identifies substantive contributions,
    explains why central claims should be accepted, describes checks and limitations,
    and responds to criticism and correction. They distinguish their contributions
    from those of advisors, collaborators, sources, and AI systems, without treating
    attribution as a substitute for mathematical understanding.
\end{enumerate}

These capacities include \term{epistemic responsibility}: an obligation to give a
defensible account of what the work establishes and why. Responsibility is necessary
for inquiry, but it is not its whole purpose. The formation of a mathematician should
also develop curiosity, imagination, technical initiative, and the ability to turn an
interesting question into a sustained investigation.

Responsibility does not require a student to reproduce every dependency personally.
Modern mathematical work already relies on established results, collaborators, and
computational systems. A student should have detailed command of their substantive
contributions and the project's central mathematical structure; identify components
on which they rely; explain the evidence that warrants that reliance; and know when
additional expertise or checking is required. The expected depth varies with the
component's importance, novelty, and risk.

This is a demanding but distributed standard. A dependency map or a verified component
can support responsible reliance, but neither can replace personal mathematical command.
A dissertation committee should be able to distinguish a student who understands and
develops a project from one who merely coordinates tools and repeats their conclusions.

\clearpage
\part{A proposed approach to training and advising}

\section{Designing learning with powerful tools}

The instructor or advisor, together with the student as their independence develops,
should identify the capacity an activity is intended to develop and choose forms of
assistance that serve that purpose. An advisor may reasonably prefer a
modest problem that develops a student's reasoning to an impressive result the student
does not understand. Equally, an ambitious AI-assisted investigation may teach more
than a long sequence of exercises whose main difficulty is avoidable technical labor.

Neither a general prohibition of AI nor its use without an educational plan is adequate.
Programs need several
learning designs: independent attempts followed by hints; study of a complete argument
followed by reconstruction and variation; collaborative exploration; and AI-open
research with subsequent explanation and evaluation. No single order is appropriate
for every topic, student, or stage.

Purposeful independent work remains important. Students need opportunities to generate
ideas, experience uncertainty, diagnose obstacles, and decide what to try without
immediately receiving a proposed next step. These opportunities should be bounded and
responsive. Indefinitely withholding explanations, hints, or feedback from an advisor,
peers, or AI can waste time, reinforce misconceptions, or make success depend too heavily
on prior preparation.

An activity can specify an initial attempt, the hints or criticism available afterward,
when a fuller solution may be consulted, and what reconstruction, adaptation, or extension
the student must then undertake. Another activity may begin with a reliable solution and
ask the student to explain or transform it. Across activities, however, reconstruction
should not become the only kind of intellectual work. Students also need occasions to
choose a direction while its prospects remain uncertain. The objective is practice in
mathematical decision-making, not a quota of failed attempts or time spent struggling.
The test is what the student learns and can subsequently do, not how closely the
sequence resembles an inherited apprenticeship.

\section{Knowing an answer as a beginning}
\label{sec:answer-first}

A reliably established answer can initiate substantial mathematical inquiry. Students
already learn by studying known theorems; AI-assisted research may create further
opportunities to begin from a result and ask why it is true, how far it extends, or what
else follows. Knowing an answer can orient these investigations without making them
routine or guaranteeing that a simple proof exists.

Consider a student studying a theorem with a technically demanding but independently
checked proof. The educational project need not be to rediscover that proof before
seeing it. Instead, the student and advisor can ask what mechanism drives the theorem.
The student uses AI to propose alternative formulations, candidate weakenings of the
hypotheses, and families of examples. These proposals remain conjectural until checked.
Some may be rejected by counterexamples; others may lead to a valid extension, a
cleaner argument in a special case, or a connection to a different method.

The student then develops one direction, explains why others failed, compares the
resulting argument with the starting proof, and identifies what has been learned.
Assessment concerns mathematical choices, technical work, explanation, and application
to new problems,
not ownership of the first answer. The activity may produce new mathematics, or it may
be valuable formation without a publishable contribution. Those outcomes should be
identified honestly.

The same design can begin with an AI-produced proof, but only after there is sufficient
evidence to treat its conclusion as established. A plausible announcement does not give
the student a known answer. If the starting proof is unverified, the project is partly
an investigation of whether the result is true, and its learning goals should say so.

Answer-first inquiry complements rather than displaces open exploration. Students also
need experience choosing questions whose answers are unknown, deciding whether to seek
a proof or a counterexample, and abandoning an unproductive formulation. Graduate
formation should include both the freedom created by a known result and the uncertainty
inherent in research.

\section{Verification, independence, and community}

Verification belongs throughout training. Early activities can test definitions against
examples or compare a generated explanation with a textbook. Later activities can
analyze limiting cases, construct counterexamples, audit dependencies, or employ
symbolic and formal methods. These are not merely inspections of completed work:
they can reveal why a result is true, how far it extends, and what questions follow.

The speed of obtaining material must be distinguished from the pace of consolidating
it. A summary of several papers is not the same as command of their ideas. Students
need opportunities to reconstruct arguments, connect concepts, teach, and apply what
they have learned. Retention, application to unfamiliar problems, and increasing initiative are more informative
than the quantity of material encountered.

Planning should itself become part of the student's education. Early on, the advisor
explains why one paper precedes another or a promising direction is postponed. AI can
suggest prerequisites and alternative routes, but suggestions must be discussed and
tested. Control should move from advisor design through co-design to student design;
meaningful student control may include extensive AI assistance.

Independence does not mean isolation. Seminars, reading groups, teaching, and contact
with several mathematicians expose students to different standards and styles. They
prevent the output of one system or preferences of one advisor from becoming the entire
intellectual environment. The aim is a student who can direct their development while
participating in a technological and human community.

\section{Roles and conditions of AI use}
\label{sec:ai-roles-conditions}

The same system can support or undermine learning depending on the role it is assigned.
Four roles are especially useful to distinguish:

\begin{enumerate}
    \item \textbf{AI as tutor.} It explains definitions, changes the level of exposition,
    produces examples, asks diagnostic questions, and offers hints.  The objective is the
    student's understanding rather than the completion of an artifact.

    \item \textbf{AI as critic.} It examines work already attempted by the student,
    points to possible gaps, requests justification, proposes tests, and compares
    approaches.  It may identify that a step is unsupported without supplying the repair.

    \item \textbf{AI as collaborator.} It proposes strategies, performs calculations,
    searches for related ideas, drafts proofs, and helps develop a project.  The student
    must track the division of labor and retain responsibility for the mathematics.

    \item \textbf{AI as object of examination.} It produces proofs, definitions,
    explanations, or plans that the student must audit, repair, compare, or reject.  The
    possibility of error is educationally central.
\end{enumerate}

These roles may be combined, but they should not be confused.  A student does not learn
proof construction merely by asking a collaborator for a proof, and a tutor should not
quietly become the unacknowledged author of assessed work.

These roles can be explored under different conditions of access to AI. A \emph{role}
describes what AI does; an \emph{access condition} describes the assistance available
during a particular activity or interval. The following are possible educational
designs, not a required classification of every activity, a sequence all students must
follow, or a prescription for particular AI systems. They may be useful in ordinary
study and research as well as assessment:

\begin{center}
\begin{tabularx}{\textwidth}{>{\bfseries}p{0.15\textwidth}Y Y}
\toprule
Condition & Possible arrangement & Intended educational purpose \\
\midrule
AI-off & A bounded interval without generative AI assistance; other available resources
are agreed separately. & Develop and demonstrate internal competence, persistence, and
live reasoning. \\
AI-limited & Selected hints, questions, criticism, or examples, intended to support the
student without supplying a complete solution. & Support the student's reasoning while
providing responsive guidance. \\
AI-open & Broad use of AI for exploration, generation, calculation, checking, and drafting,
with appropriate disclosure. & Develop responsible collaboration with AI and pursue
authentic or ambitious research. \\
\bottomrule
\end{tabularx}
\end{center}

Which arrangements are useful depends on the mathematical objective, the student's
development, and the setting. An advisor and student might agree on a bounded period
of independent exploration, selected assistance while the student develops an approach,
or an investigation that begins from a complete argument. These are possibilities to
adapt, not a fixed sequence or a new system of surveillance.

Their practical limits matter. An agreement about private study is not proof that a
restriction was observed, and asking a model to give hints does not guarantee that it
will withhold a solution. Where restricted assistance is important to an assessment,
the program needs a practicable arrangement and appropriate accommodations; it should
not claim evidence that the arrangement cannot provide. The aim is to support learning
and obtain proportionate evidence of development, not to make every activity fit a
category. Broad AI use still requires mathematical responsibility and attention to
privacy and equitable access.

Section~\ref{sec:assessing-capacities} discusses evidence of formation, and
Appendices~\ref{app:protected-work} and~\ref{app:ai-open} offer illustrative learning
cycles. Their usefulness should be evaluated rather than assumed.

\section{An integrated training model}

Graduate formation should not be organized as ``learn first, verify later, research
last.''  Several strands should run through the program at increasing levels of
complexity and student control.

\subsection{Knowledge, technique, and mathematical orientation}

Students still need substantial internal mathematical knowledge.  Background knowledge
determines which questions they can formulate, which analogies they notice, which outputs
appear suspicious, and which explanations they can integrate.  Foundational courses
should remain demanding while distinguishing three legitimate educational goals:

\begin{itemize}
    \item \term{working fluency}: definitions, examples, core theorems, and techniques
    that should be available without extensive external assistance;
    \item \term{navigational knowledge}: awareness of areas, methods, mathematicians with relevant expertise, and
    references sufficient to know what to seek and where to seek it;
    \item \term{tool-mediated capability}: computations, searches, formalizations, or
    specialized procedures that may be performed with external systems, but whose
    assumptions and outputs the student can evaluate.
\end{itemize}

Fields will draw these boundaries differently.  The important point is that they be
chosen for educational and mathematical reasons rather than inherited from current tool
limitations.

Proof, computation, and example-building remain central.  Constructing a proof teaches
more than the final sequence of implications; it develops control of definitions,
assumptions, and possible approaches.  Examples and counterexamples reveal how objects
behave and resist plausible overgeneralization.  Students should routinely ask which
hypotheses are used, what occurs in limiting or degenerate cases, and whether scaling,
dimensional, numerical, or other domain-appropriate checks support a claim.

\subsection{Guided study and research apprenticeship}
\label{sec:guided-study}

The advisor introduces a research area through selected papers, examples, and questions.
AI can help map prerequisites, offer explanations at different levels, and generate
practice material.  The advisor remains responsible for the intellectual design and for
explaining the choices.  The student presents material independently in meetings,
answers questions, reconstructs arguments, and identifies what remains uncertain.

When a research problem is introduced, the advisor, student, and AI may decompose it into
subproblems.  The decomposition should be visible and contestable: the student should
understand why each part matters and gradually learn to propose the decomposition rather
than merely follow it.  One cycle moves from an independent attempt to bounded assistance, testing,
explanation, and selection of the next question. Another starts from an established
result and develops a new direction through AI-open exploration. Appendices~\ref{app:protected-work}
and~\ref{app:ai-open} illustrate these designs; neither should become a compulsory
routine.

The advisor adjusts support based on observed learning. 
A brief attempt made only to satisfy a requirement before requesting a full AI-generated answer does not demonstrate meaningful learning. 
Conversely, extensive AI use does not by itself indicate a lack of mathematical development. 
What matters is the student's substantive participation and the capacities that participation develops. 
As students mature, they should gain greater freedom to choose
their tools and methods, together with greater responsibility for the mathematical
decisions that follow.

\subsection{Progression toward research leadership}

The training model describes a progression in responsibility rather than a fixed
calendar.  At entry, students establish foundational fluency and learn the program's
standards for AI use, verification, and attribution.  As they enter a research area, they
read and present key work, construct examples, and participate in designing their
learning plans.  In guided research they may explore different arrangements for assistance, develop
domain-specific forms of verification, and keep concise records of important decisions
and contributions.

In the final stage, the student increasingly selects questions, designs the project,
chooses the division of labor, and decides what evidence is sufficient.  Success should
not be defined as producing something that ``AI alone could not have produced,'' an
unstable and unmeasurable counterfactual.  The relevant standard is that the student's
sustained participation adds identifiable intellectual value and that the student can
take responsibility for the final work.

Mathematical maturity also includes communication and teaching.  Explaining one result to
a specialist, a mathematician outside the area, and an advanced undergraduate reveals
whether the student sees the underlying structure or only reproduces notation.  Teaching
and exposition are therefore not secondary presentation skills; they provide evidence of
understanding and prepare students to participate in a community of knowledge.

\section{The advisor in the new model}
\label{sec:advisor-model}

The advisor helps students develop mathematical imagination, technical command, taste,
and research independence. This includes selecting reliable sources, discussing
worthwhile problems, making the reasons for a research choice visible, and recognizing
when a student's unexpected direction deserves support. The advisor also distinguishes
genuine understanding from fluent output and models how to use AI critically and
creatively. Establishing that the student genuinely understands an argument is not
enough: the advisor must also consider whether the student is developing the ability
to initiate and direct mathematical work.

The frequency of research meetings need not change. What needs reconsideration is the
intellectual activity between meetings and the evidence discussed within them. Alongside
completed arguments, meetings should sometimes address an unfinished question: what
does the student think is worth trying, what alternatives are plausible, and what evidence
would justify changing course? Discussing choices before their outcomes are known
provides different evidence from a polished retrospective explanation.

For a selected task, student and advisor can identify the capacity to be developed and
agree on a suitable division of work. The student might formulate candidate lemmas
before asking AI for criticism, compare assisted approaches and choose one to pursue,
or start from a supplied proof and design an investigation of an unresolved extension.
A brief note of a consequential choice and its reasons can support the discussion;
routine progress need not generate a new reporting burden. The next meeting considers
both the mathematics and whether the activity exercised the intended capacity.

These practices do not make the advisor a monitor of private study or require withholding
help until the student fails. They make AI use part of the deliberate design of the
apprenticeship. If existing advising already provides these opportunities and credible
evidence of development, it may need little structural change. What cannot be assumed
is that project progress plus genuine understanding establishes that all is well.

AI may reduce some routine demands and increase others. Generating practice material
or clarifying prerequisites may save time, while reviewing sophisticated generated
arguments may require more. The balance should be observed rather than assumed.
Whatever the division of labor, the student needs substantive mentoring, professional
advocacy, and access to mathematical perspectives beyond one system or one advisor.

An AI-generated study plan may explain why a prerequisite should be studied first,
why one paper is relevant, or why a particular problem is a useful next step. Those
explanations can support learning, but the advisor and student should still discuss
whether they fit the student's needs and research goals. The aim is not to reserve
planning for humans, but to develop the student's ability to evaluate, revise, and
eventually direct a plan.

Preparing faculty members to advise and teach in this environment should involve actual
mathematics: exploring AI-assisted approaches, comparing proofs, designing activities
with different kinds of assistance, and discussing how to assess understanding fairly. It should also address attribution, confidential material, unequal preparation,
and the distinction between educational observation and surveillance. Learning records
should support feedback, not demand continuous monitoring of a student's private work.

The advisor must protect the educational purpose of the relationship. Verification is
part of mathematical inquiry and may itself be creative research, but assigning students
routine checking of an advisor's generated output is not an adequate apprenticeship.
Their work should progressively include question formulation, method development, and
ownership of a research direction.

\clearpage
\part{Assessment and program redesign}

\section{Assessing capacities rather than artifacts alone}
\label{sec:assessing-capacities}

Research results remain valuable.  They can show ambition, persistence, and the ability
to bring a project to completion.  They can no longer serve as the primary proxy for
student development.  Programs need direct evidence of competence, judgment,
independence, and responsibility.

\begin{center}
\begin{tabularx}{\textwidth}{>{\bfseries}p{0.18\textwidth}Y Y}
\toprule
Capacity & Illustrative evidence & What the evidence helps reveal \\
\midrule
Competence & Proof construction and reconstruction, examples, computation, development
of a new approach, application to unfamiliar problems & Whether knowledge and technique support productive reasoning \\
Judgment & Flaw detection, comparison of approaches, confidence with reasons, evaluation
of significance and naturalness & Whether the student can distinguish plausibility,
correctness, and mathematical value \\
Independence & Student-formulated questions, research plans, creative use of assistance,
choice and revision of methods & Whether the student can initiate and develop inquiry \\
Responsibility & Contribution statement, checks of central claims, account of uncertainty,
oral defense of decisions & Whether the student can understand and stand behind work
presented under their name \\
\bottomrule
\end{tabularx}
\end{center}

No single instrument can measure all four capacities.  Assessment should combine work
produced over time with direct demonstrations under clearly described conditions.
The possibilities in Section~\ref{sec:ai-roles-conditions} can serve different purposes
where they are practical: AI-off work may show foundational fluency and live reasoning;
AI-limited work may show how the student uses guidance; AI-open work may show responsible
collaboration in an authentic research environment. No single activity must use all
three arrangements.

Take-home work remains useful for assessing exposition, synthesis, revision, and tool use,
but it becomes weak evidence of unaided competence when considered alone.  A short oral
follow-up, whether part of an examination or an ordinary meeting with the advisor, can
ask the student to explain a key step, adapt the solution to a changed hypothesis,
identify where assistance entered, or reconstruct the argument without the submitted
text.  Oral examination should not become the sole safeguard: quick verbal
performance can reflect language background, disability, anxiety, and style as well as
understanding. Observations across regular advising meetings and seminars, written work,
and sustained projects provide complementary evidence. Assessment should include
appropriate accommodations rather than depend on one style of performance.

Even an excellent explanation or reconstruction cannot by itself establish research
initiative. Evidence should also include work on questions before a route has been
supplied: a proposed formulation, reasons for trying an approach, a revision prompted
by an obstacle, or the design of a next investigation. Some of this work can be AI-open;
the issue is the student's substantive participation in making decisions, not the
absence of help. Discussing what the student intends to try before the outcome is known,
and observing how they approach new tasks over time, provides evidence that a later
account of a completed project cannot supply by itself. No decision log or single
successful demonstration establishes independence on its own.

\section{Exercises in judgment and creative inquiry}

Several assessment practices are especially suited to an environment filled with
plausible AI-generated mathematics.  They can be used formatively before they carry
high-stakes consequences.

In a \term{flaw-detection exercise}, an instructor begins with correct arguments and
introduces controlled errors: a dropped hypothesis, reversed quantifiers, an unjustified
extension to a boundary case, an incorrect dependency of a constant, an unjustified exchange of limits, an
unhandled case, or a vacuous conclusion.  Students locate and classify the issue, explain
why it matters, state their confidence, and repair it if possible.  Correct arguments
must be included as negative controls; otherwise students are rewarded for inventing
objections to every proof.  Results should distinguish detection from false positives and
should be examined by error category.

Confidence itself can be trained.  For selected judgments, students report a coarse
confidence level together with the evidence behind it.  Across many tasks, the program
can help students see whether they are systematically overconfident, or whether they
identify genuine problems but distrust correct arguments as well.  The objective is not
to attach decorative percentages to statements, but to develop calibrated self-knowledge.

Other useful exercises ask students to remove a hypothesis and construct a counterexample,
analyze a boundary case, compare several correct approaches, rank research questions by
significance and defend the criteria, or explain the same result at several levels.  Such
activities reveal forms of understanding and judgment that a polished final artifact may
conceal.

Creative tasks should be equally visible. Students can formulate a conjecture from
examples, design an informative computational experiment, propose two approaches to a
question, or develop a connection suggested by AI. Assessment should consider the
quality of the question, mathematical reasoning, use of evidence, and productive
revision, without rewarding confident novelty claims over careful investigation.

\section{Milestones, examinations, and the dissertation}
\label{sec:milestones}

Qualifying examinations should state what they certify. AI-off components can test
foundational fluency, AI-open components can test synthesis and productive tool use,
and oral or written follow-ups can test understanding and adaptation. No single mode
should carry the entire burden of establishing mathematical maturity.

Candidacy should demonstrate the ability to enter a field, understand and examine a
substantial argument, formulate a worthwhile question, and design a plausible
investigation. A proposal should explain the role of AI and other resources, what
uncertainties remain, and how the student will respond to evidence that the original
direction is unpromising.

\subsection{Originality and the student's contribution}

The dissertation should make a substantial contribution to mathematics and demonstrate
the student's capacity for independent research. These are related but distinct
requirements. An important AI-assisted result does not automatically establish the
student's maturity; excellent understanding of existing work does not by itself
constitute an original research contribution.

Originality need not be confined to obtaining the first proof of a new theorem. A new
method, a significant generalization, an illuminating alternative proof, a conceptual
connection, or a mathematically substantive formalization can advance a field.
These are possibilities for evaluation, not automatic qualifying categories.
Routine formalization, competent exposition, or a useful compilation may be valuable
without meeting the standard for a doctoral research contribution.

The committee should identify what is new, why it matters, and how the student
contributed intellectually to its development. Contributions may include problem
formulation, mathematical design, technical arguments, interpretation, or a decisive
revision of an approach. They must be substantial and demonstrable across the project,
not merely a claim to having initiated an AI interaction or assembled its outputs.

The standard is not that the student has done something AI could never do. It is that
the work advances mathematics and the student's sustained participation demonstrates
original, independent, and responsible inquiry. As elsewhere in collaborative research,
credit for a shared achievement and evidence of an individual's formation require
careful, explicit judgment.

\subsection{Evidence, defense, and proportionate records}

The dissertation should be accompanied by a concise responsibility and contribution
statement addressing main claims and important decisions rather than every sentence.
It identifies substantive human and AI contributions, relevant sources, checks,
dependencies, unresolved limitations, and enough information to examine important
computational or formal components.

The defense tests whether the student can explain, modify, criticize, and contextualize
the work. A committee might change a hypothesis, request a special case, ask how a new
direction was chosen, or examine a component on which the project relies. It should
also test creative understanding: what further question does the work make possible,
and what would a serious attempt to answer it involve? These demonstrations should be
considered alongside evidence of initiative and changing responsibility during the
research, not treated as a substitute for it.

Records should remain proportionate. Useful details may include the name and release
of an AI model, the versions of computational software or proof assistants used, relevant
settings, and selected interactions that materially shaped the work. An archive of
every query is neither necessary nor sufficient. Repeating a query may not reproduce
the same response, and an earlier model may no longer be available. What must remain
open to examination is the mathematics itself: the proof or code actually used, its
assumptions and dependencies, and the checks that support its conclusions. Where
supporting material is confidential, authorized examiners should have a secure way to
review what they need without requiring that the material be made public.

\section{Curriculum, policy, and community}

Core mathematical education should not disappear because information is retrievable.
Courses should identify the concepts and techniques students need internally, emphasize
connections and the application of learned ideas to unfamiliar problems, and choose
forms of AI assistance for specific educational purposes.  Breadth may become more important, not less, when students must evaluate
connections and generated analogies across areas.

Verification and judgment should form a thread through the program rather than a single
new course.  A first-year class can use elementary flaws and counterexamples; an area
seminar can examine research-level arguments; an advanced group can referee AI-assisted
drafts or use formal and computational checks where relevant.  Supervised peer-review
seminars can distribute practice economically, but faculty must set standards, audit the
discussion, and correct shared errors.

Programs also need task-specific policies.  A single rule cannot govern homework,
qualifying examinations, reading courses, and open research equally.  Policies should
state the permitted role of AI, disclosure expectations, data restrictions,
responsibility for errors, and consequences of misrepresentation.  Faculty should
disclose their own substantive AI use when collaborating with students; responsibility
cannot be a one-way demand.  AI systems can respond to criticism and help revise an argument, but those
responses do not discharge the professional accountability of the human authors or
committees who endorse it. Institutions must assign responsibility for investigating
errors, correcting the record, and protecting affected students.

Equal access to an interface does not guarantee educational equity.  AI can provide
responsive explanation to students with uneven prior opportunities, but those with weaker
preparation may also be least able to detect plausible errors.  Paid tools, computing
access, disability, language, and variation in advisor engagement can divide a cohort.
Programs should provide a supported baseline of AI and computing access, early
verification training, multiple forms of assessment, human support for foundational gaps,
and routes for challenging automated judgments. Access to AI and computing resources
means sufficient usage allowances, time on shared systems, and technical support, not
merely an account. Shared university or consortium provision can support ordinary work,
with an additional allocation process for projects requiring substantially more
computing time, memory, or specialized technical support. Assessment must not silently depend on a student's ability to purchase
more powerful assistance.

Finally, programs should resist the apparent efficiency of replacing seminars and human
discussion with individualized tutoring.  Students need to see experts disagree, watch
questions being formulated, learn disciplinary norms, and encounter several styles of
reasoning.  Cross-area seminars, rotations, reading groups, teaching, and informal
conversation remain central to the formation of judgment.

\section{Duration, scale, and professional preparation}

The case for automatically shortening the PhD assumes that its product is research
output.  If the purpose is mature mathematical inquiry, the binding constraint may
be human consolidation rather than generation speed.  Reading, forgetting, returning,
teaching, failing, and revising are not obviously compressible because an AI system can draft
a proof quickly.

Programs should become more competency-based and less dependent on a uniform chronology.
A well-prepared student who demonstrates the required capacities should not be held only
to satisfy a customary number of years; a student should not be rushed because AI made a
dissertation draft arrive early.  Competency-based progression requires funding and clear
milestones so that ``maturation takes time'' does not justify indefinite enrollment.

Program scale is a related but separate question. Preserving a pipeline of trained
mathematicians does not determine a particular cohort size, nor does AI provide an
automatic argument for reducing it. Departments must consider advising capacity,
educational quality, career outcomes, and the full costs of student support, mentoring,
teaching provision, and computational infrastructure. Maintaining or expanding a cohort
may be justified by the institution's educational and research mission, the need to
sustain expertise, and credible opportunities for graduates, provided the necessary
support can be funded. Where support is inadequate, the alternatives include additional
investment, shared resources, and redesigned advising as well as changes in admissions.
Neither preserving existing numbers nor shrinking them substitutes for examining what
the program can responsibly provide and what purposes it should serve.

The changing profession need not mean only declining opportunity.  As AI produces more
mathematical arguments, models, computations, and conjectures, new needs may arise for
people who can understand, organize, evaluate, and responsibly apply that output.
Alongside existing careers in university research and teaching, industry, and public
institutions, mathematically trained graduates may contribute to emerging or expanded
forms of work. These include translating AI-generated ideas into reliable scientific,
engineering, and computational methods; providing mathematical oversight in
consequential applications; supporting interdisciplinary modeling; developing
AI-assisted mathematical education; and leading long-term human--AI projects requiring
problem selection, verification, synthesis, and interpretation. Some of these activities
already belong to mathematical careers; AI may change their scale, content, or
institutional setting rather than create entirely new occupations.

The precise shape and scale of these roles cannot yet be known.  Programs should not
promise them as guaranteed career paths, but neither should they assume that advanced
mathematical training will lose its value.  Students need honest discussion of risks and
possibilities together with experience in rigorous reasoning, computation, formalization
where appropriate, reproducible work, technical communication, collaboration, project
design, and learning unfamiliar domains.  Career preparation should be integrated into
authentic projects, teaching, interdisciplinary collaboration, or internships rather than
added as a generic workshop at the end of the degree.

\clearpage
\part{Institutional responsibility and implementation}

\section{Recognizing the work of advising}
\label{sec:advising-work}

Faculty performance is generally evaluated through research, teaching, and service, yet
graduate advising does not fit neatly within any one category.  In many departments it is
counted as service, treated as an additional educational responsibility on top of the
formal teaching load, valued indirectly as evidence of an active research program, or
recognized through some combination of these approaches.

These conventions may become inadequate as the work of advising changes. Some
activities may become less demanding, others more intensive, and AI-assisted students
may contribute to research in new ways. Workload policy should reflect actual mentoring
rather than presume either an automatic saving or an inevitable increase. Where
the intensive work of developing students' mathematical capacities is necessary but
unrecognized, faculty face an incentive to withdraw
from advising or provide only nominal supervision.  Departments and universities should reconsider how advising is credited in
workloads, annual reviews, salary decisions, and promotion.

One possible mechanism is a \term{graduate mentoring allocation}.  A department defines
the activities that constitute intensive formation---an individual learning plan,
regular substantive meetings, opportunities to develop research initiative, assessment
of understanding and mathematical choices, feedback on verification, and an annual
development review---and assigns documented workload credit.
Accumulated credit may produce a course release, summer support, reduced committee work,
or another locally meaningful adjustment.  Co-advisors divide rather than duplicate the
credit, and promotion review treats the mentoring record as educational contribution.

This mechanism has real costs.  A course release requires replacement teaching, fewer
offerings, larger classes, or additional budget.  The cost should be visible rather than
transferred silently to individual faculty.  Team advising, shared area seminars, and
verification workshops can reduce duplication and expose students to several standards,
but they also require coordination. If a department cannot support responsible formation
for its current cohort, it should examine additional investment, shared advising, and
changes in program design alongside future admission numbers and program promises.
Existing commitments to students should be protected; reducing admissions is not the
only response to inadequate support.

\section{A staged departmental implementation}
\label{sec:staged-implementation}

A department need not redesign the entire PhD at once, but staged implementation should
organize prompt action rather than postpone it. Departments can begin discussing AI use
in advising, clarifying responsibilities, supporting equitable access, and examining
evidence of student development while more specific educational practices are tested.
The stages below can overlap; they need not become a sequence of lengthy approvals.

A first stage can establish a common language for the four capacities in
Section~\ref{sec:capacities} and the four AI roles and possible access conditions in
Section~\ref{sec:ai-roles-conditions}. Two existing graduate courses might each try a
purposeful independent attempt, a flaw-detection exercise, and an AI-assisted creative
investigation, choosing practical arrangements suited to their goals. One activity
could begin with a reliably established answer and assess the student's ability to
explain, vary, or extend it. Selected take-home work can be paired with brief oral
follow-up. A departmental working session for advisors and instructors could examine
AI-generated arguments from their own fields, compare their educational uses, and
discuss how to assess the student's contribution and understanding.

A second stage can adapt promising practices for milestones: a qualifying component
that draws on more than one practical arrangement for assistance, a candidacy exercise
in auditing and research planning, and a contribution statement attached to advanced
projects.  The department can establish a supervised
student-refereeing seminar or incorporate peer review into an existing seminar.  At the
same time, it should review tool access, privacy rules, accommodations, advisor training,
and the workload credit attached to early-stage formation.

A third stage can reconsider the curriculum and program as a whole in light of what the
earlier stages reveal.  This may include changes to course requirements, breadth,
dissertation expectations, degree duration, cohort size, and professional pathways.  The
distinction is between acting promptly and committing prematurely. Departments should
make reversible improvements now, set explicit dates for reviewing their effects, and
revise them as evidence accumulates. Changes with lasting consequences---especially
those affecting degree requirements, admission commitments, or student funding---need
broader deliberation and protection for current students. The aim is timely action
with accountable revision, not an indefinite pilot or a demand for certainty before
anything changes.

\section{Learning from a pilot}
\label{sec:pilot}

The model contains empirical assumptions: that purposeful independent work supports
retention and the use of learned ideas in new settings, that verification practice
improves judgment, and that guided
AI collaboration can develop both mathematical initiative and effective tool use.
These assumptions should be examined rather than embedded permanently in policy.

A first pilot should evaluate components, not compare tiny graduating cohorts trained
under entirely different regimes. In courses or reading groups, a pilot might compare
matched tasks under selected AI-off, AI-limited, or AI-open arrangements where those
arrangements can be implemented credibly. Students can encounter the arrangements in
different orders to reduce the influence of practice or sequence. The pilot need not
include all three. Another comparison could examine an initial attempt versus beginning
with an established solution. Tasks should be comparable in difficulty; all students
should have access to the learning opportunities and appropriate accommodations.

Evaluation should distinguish three outcomes:
\begin{enumerate}
    \item \textbf{Immediate mathematical work:} correctness, quality of explanation,
    quality of questions or approaches, useful revisions, and time required.
    \item \textbf{Internal formation:} reconstruction after a delay and application
    to unfamiliar problems without the original output, together with formulation and
    choice of approaches on a new question before a route is supplied.
    \item \textbf{Assisted research capability:} later work on an unfamiliar problem
    with AI available, assessing how students choose assistance, test suggestions,
    develop ideas, and revise direction.
\end{enumerate}

The third outcome matters because a program devoted to human--AI inquiry should not
evaluate success exclusively through unaided performance. Conversely, a strong
assisted result cannot alone establish what a student has internalized. The two kinds
of evidence are complementary.

Judgment tasks should track false positives as well as error detection, with results
examined by error category. Creative tasks need explicit criteria and examples of
acceptable reasoning, so that evaluation does not merely reward an assessor's preferred
style. Where feasible, assessors should review work without knowing its instructional
condition, and disagreements should inform refinement of the criteria.

The pilot should also record student burden, advisor and staff time, computational
resources, and differences associated with preparation or access needs. A practice
that works only with exceptional supervision may need redesign before broad adoption.
Research participation and data use require appropriate institutional review and
consent procedures; ordinary educational support should not depend on agreeing to
additional research data collection. Small samples and changing systems limit what
can be concluded, so the conditions and relevant tools should be documented.

Before the pilot, the department should state how findings would affect its plans.
If independent attempts improve later application to new problems at reasonable cost, they can be
expanded where effective. If AI-open work produces equal or better formation, restrictions
should be narrowed. If one design improves assisted investigation while weakening
foundational command, the program should investigate how to combine its benefits with
targeted practice. If flaw-detection work raises both valid detections and false
accusations, calibration and evidentiary standards need attention.

Such a pilot cannot establish the value of an entire PhD or predict careers. It can
identify promising practices, expose attractive assumptions that fail, and reveal the
resources needed for responsible implementation.

\clearpage
\section{Conclusion}

AI makes it necessary to distinguish more carefully between producing mathematical
work, understanding it, and developing the capacity to initiate and direct inquiry.
A valuable result and a genuinely understood proof can coexist with insufficient
practice in formulating questions and choosing approaches under uncertainty. But
assistance can also expand the student's opportunities for exploration, explanation,
and discovery. Graduate education should deliberately shape that participation rather
than infer complete formation from progress and understanding alone.

The central objective is original, independent, and responsible inquiry. It requires
substantive knowledge and technical skill, the imagination to formulate questions and
develop approaches, judgment about evidence and significance, and the capacity to
collaborate without becoming a passive recipient of conclusions. Mathematical maturity
is demonstrated in what a student can initiate, understand, build, explain, revise,
and responsibly contribute.

A redesigned apprenticeship should therefore combine purposeful independent work,
learning from established answers, creative AI-open projects, and continuing practice
in verification. It should distinguish mathematical contribution from evidence of
formation, make doctoral originality explicit, and permit justified reliance on others
without surrendering personal command. Institutions must provide equitable support,
recognize the actual work of advising, and assess whether their educational designs
succeed.

The aim is not to preserve difficulty for its own sake or to locate a shrinking set of
tasks beyond the reach of AI. It is to educate people who can participate deeply in
mathematics as its possibilities expand. The methods should remain open to revision;
the ambition is a richer formation that enables mathematicians to ask more, understand
more, and contribute more through their own work and their collaboration with powerful
tools.

\section*{AI assistance disclosure}
\addcontentsline{toc}{section}{AI assistance disclosure}
This article was developed with assistance from ChatGPT/Codex and Claude. These tools supported the exploration and refinement of arguments, organization of the manuscript, and revision of the text. The author made the final decisions about the positions advanced and is responsible for the article’s content.
\clearpage
\appendix
\small
\setlength{\parskip}{0.35em}
\setlist[enumerate]{leftmargin=1.65em,topsep=2pt,itemsep=2pt,parsep=0pt}

\section{A sample protected-work learning cycle}
\label{app:protected-work}

The following cycle combines protected work, AI assistance, verification, and advisor
feedback for a student entering a research area. Its timing and difficulty should be
adjusted to the student and material; it is one design, not a required schedule.

\begin{enumerate}
    \item \textbf{Orientation.} The student reads the advisor's statement of the activity's
    objective and writes what they already know, what appears unfamiliar, and why the
    topic matters to the larger plan.
    \item \textbf{Protected reading and attempt.} The student reads selected primary or
    textbook material and reconstructs a central example or proof without AI 
    assistance for a specified interval.
    \item \textbf{Tutor or critic interaction.} AI asks questions, supplies alternative
    examples, or gives bounded hints about points the student has identified.
    \item \textbf{Independent reconstruction.} With the interaction closed, the student
    rewrites the argument, performs a related computation, or applies the idea to a new problem.
    \item \textbf{Audit.} The student compares the reconstruction with trusted sources,
    tests boundary cases, and records what remains uncertain.  A short generated argument
    may be included specifically for criticism.
    \item \textbf{Advisor or group meeting.} The student explains the mathematics
    without a generated script and discusses a choice or obstacle, including one
    unresolved question and what they would try next. The advisor assesses understanding
    and initiative, then adjusts the plan.
    \item \textbf{Increasing student control of planning.} The student proposes the next objective and the
    appropriate access condition.  Early in the program the advisor may substantially
    revise the proposal; later, the student should control it.
\end{enumerate}

\section{A sample AI-open inquiry cycle}
\label{app:ai-open}

This cycle complements the protected-work design. It can begin with a trusted theorem
or with an open question; the mathematical status of the starting material must be clear.

\begin{enumerate}
    \item \textbf{Choose a direction.} The student states a question and why it matters.
    If beginning from a known theorem, they identify a possible extension, explanation,
    alternative proof, or consequence worth investigating.
    \item \textbf{Explore with assistance.} The student uses AI, literature, computation,
    or collaborators to develop examples and candidate approaches. Generated claims
    are recorded as proposals, not established facts.
    \item \textbf{Make mathematical choices.} Before the outcome is known, the student
    selects a direction and gives reasons for trying it. They do substantive work,
    revising the plan as evidence develops. This may include proof construction,
    counterexample design, a new formulation, or a decisive comparison.
    \item \textbf{Establish and explain.} The student checks what the investigation
    supports, identifies dependencies, and separates established conclusions from
    unresolved questions.
    \item \textbf{Discuss and extend.} In a meeting, the student explains the main
    decisions, responds to a changed assumption or new objection, and proposes a next
    question. The advisor assesses both personal command and effective collaboration.
\end{enumerate}

A concise record of consequential decisions is more useful here than a complete
interaction log. A project that yields no new theorem can still show significant
formation; any claim to new mathematics must be evaluated separately.

\clearpage
\section{A sample flaw-detection assessment}

An assessment set contains correct and modified arguments matched by topic and approximate
difficulty.  For each argument the student submits:

\begin{enumerate}
    \item a verdict: no identified problem, minor issue, repairable gap, serious problem,
    or unable to determine;
    \item the exact location and category of any issue;
    \item an explanation or counterexample;
    \item a possible repair, if one is available;
    \item a coarse confidence level and the evidence supporting it;
    \item the additional information or check that would most reduce uncertainty.
\end{enumerate}

Possible injected categories include dropped or weakened hypotheses, swapped quantifiers,
incorrect constant dependencies, invalid boundary cases or scaling, unjustified exchanges of
limits or derivatives, unhandled cases, circular dependencies, invalid citation use,
vacuous conclusions, and irrelevant arguments that do not establish the stated claim.

The scoring report should separate localization, diagnosis, severity, repair, confidence,
and false-positive rate.  The instrument should be used formatively before it is used for
high-stakes evaluation.

\section{A responsibility and contribution statement}

For a dissertation chapter, paper, or major project, the student can address the
following questions:

\begin{enumerate}
    \item \textbf{Purpose and significance:} What question does the work address, and why
    is it worth addressing?
    \item \textbf{Human contributions:} Which central ideas, decisions, proofs,
    computations, explanations, and revisions came from the student, advisor, or other
    collaborators?
    \item \textbf{AI contributions:} Which systems were used, in which roles, and for
    which substantive components?
    \item \textbf{Sources and novelty:} How were related results located, and how was the
    relationship to existing work assessed?
    \item \textbf{Verification and reliance:} For each main result, what checks were
    performed, and what remains uncertain? Which dependencies does the student
    understand in detail, which are relied upon, and what warrants that reliance?
    \item \textbf{Rejected paths:} What important generated suggestions were rejected or
    substantially changed, and why?
    \item \textbf{Limitations:} Which assumptions, boundary cases, dependencies, or open
    gaps should a reader know about?
    \item \textbf{Reproducibility and privacy:} What supporting material can be shared,
    and what cannot be disclosed because of privacy, licensing, confidentiality, or size?
\end{enumerate}

\end{document}